\documentclass[a4paper,11pt]{article}
\usepackage[utf8x]{inputenc}
\usepackage[english]{babel} 
\usepackage{ucs}
\usepackage{textcomp}
\usepackage{dsfont}
\usepackage{natbib}
\usepackage[T1]{fontenc}
\usepackage{amsmath, amsfonts, amsthm, amssymb}
\usepackage{fullpage}
\usepackage{float} 
\usepackage{array}
\usepackage{hyperref} 
\usepackage{siunitx}
\usepackage{graphicx}
\usepackage{url}
\usepackage{subfig}

\newcommand{\itemr}{\item[$\rightarrow$]}

\newcommand{\R}{\mathbb{R}}

\newcommand{\dx}{\textrm{d}}
\newcommand{\bx}{\mathbf{x}}
\newcommand{\be}{\mathbf{e}}
\newcommand{\bv}{\mathbf{v}}
\newcommand{\bu}{\mathbf{u}}
\newcommand{\bn}{\mathbf{n}}
\newcommand{\Div}{\mathrm{div}\;}
\newcommand{\co}[3]{\tilde{c}^{#1}_{(#2,#3)}}

\title{Numerical solutions of a 2D fluid problem coupled to a nonlinear non-local reaction-advection-diffusion problem for cell crawling migration in a discoidal domain}
\author{Christ\`ele Etchegaray\thanks{Institut de Math\'{e}matiques de Toulouse, Universit\'{e} Paul Sabatier, 118, route de Narbonne
F-31062 Toulouse Cedex 9 ({\tt christele.etchegaray@math.univ-toulouse.fr})}  \and Nicolas Meunier\thanks{MAP5, CNRS UMR 8145, Universit\'{e} Paris Descartes, 45 rue des Saints  P\`{e}res
75006 Paris,
France. ({\tt nicolas.meunier@parisdescartes.fr})}}
\date{}

\begin{document}

\maketitle

\section*{Abstract}
In this work, we present a numerical scheme for the approximate solutions of a 2D crawling cell migration problem. The model, defined on a non-deformable discoidal domain, consists in a Darcy fluid problem coupled with a Poisson problem and a reaction-advection-diffusion problem. Moreover, the advection velocity depends on boundary values, making the problem nonlinear and non local. \par
For a discoidal domain, numerical solutions can be obtained using the finite volume method on the polar formulation of the model. Simulations show that different migration behaviours can be captured.

\begin{description}
\item[Keywords:] 2D cell migration, Darcy fluid, nonlinear non-local reaction-advection-diffusion problem, finite volume method. 
\end{description}

\section{Introduction}
Cell migration ensures fundamental functions in the body (embryogenesis, immune system), but is also involved in the development of pathologies such as tumor metastasis, arising large research efforts \citep{Bravo-Cordero2012Directed-cell-i}. However, the responsible intracellular mechanisms involve multiscale interaction in time and space, so that modelling cell migration is challenging and produces interesting problems to study. \par 
For 2D cells crawling on a surface, the motion is friction-based and results from the activity of the so-called \textbf{actin cytoskeleton}, that is a dynamics mesh of actin filaments. They are polar: they polymerize at one end and depolymerize at the other end, under the molecular regulation of signalling loops. Overall, the actin mesh can be approximated by an active fluid in our setting \citep{Kruse2005Generic_theory,Joanny2009Active}.

The mesh grows preferentially at the cell membrane, and shrinks inside the cell body, generating inward actin flows from the membrane. Its mechanical connection to the adhesive substrate generates the friction forces responsible for the cell's displacement. 

Modelling this process is a difficult task, because of the large time and space scales, and also because of the large number of effectors of the dynamics. Following key physical ideas of \cite{Blanch-Mercader2013Spontaneous-Mot,Maiuri2015Actin-flows-med}, we proposed in \cite{etchegaray:tel-01533458,etchegaray2017analysis,JMB2D} a minimal multiscale model for 2D crawling migration, that we recall now. \par 

The cell domain is a non deformable disc $\Omega$, and the problem is formulated in the cell's frame of reference (the domain does not move). First, the actin cytoskeleton is approximated by a Darcy fluid \citep{Blanch-Mercader2013Spontaneous-Mot}, and its dynamics in a crawling situation is modelled by a Poisson problem on the fluid pressure \citep{etchegaray2017analysis}. More precisely, for $\bu: \R_+ \times \Omega \rightarrow \R^2$ the fluid velocity, $p: \R_+ \times \Omega \rightarrow \R$ its pressure, we have 

\begin{equation}\label{Poisson}
\begin{cases}
            -\Delta p(t,\bx)  =  -k_d  & \text{ in } \Omega\, ,  \\
			p(t,\bx) = k_p(t,\bx)  &  \text{ on } \partial \Omega\, , 
\end{cases}
\end{equation}
for $k_d$ the actin depolymerization rate, and $k_p: \R_+ \times \Omega \rightarrow \R$ the polymerization function rate at the boundary. The fluid velocity writes
\begin{equation}\label{Darcy}
\bu(t,\bx) = - \nabla p(t,\bx)  - \bv(t)      \quad \text{ on } \Omega\, , 
\end{equation}
with for $\gamma \in \R_+$ the domain velocity
\begin{equation}\label{DomVeloc}
 \bv(t) = \gamma \int_{\Omega} \nabla p(t,\bx) \,  \dx \bx = \gamma \int_{\partial \Omega} p(t,\bx) \bn_{\bx}\, \dx \bx \,,
\end{equation}
for $\bn_{\bx}$ the unit normal vector to $\partial\Omega$ at point $\bx \in \partial \Omega_0$.

These equations show that the actin polymerization at the boundary and depolymerization inside the domain may drive a pressure gradient leading to the cell motion. Note that the dynamics therefore arises from the activity at the boundary, and that the cell velocity is nonlocal. 

Now, the interaction with the molecular scale consists in studying the dynamics of a molecular inhibitor to polymerization. The molecules diffuse freely inside the cell in a inactive form. They may bind actin filaments and be carryied by their flow. Moreover, if activated at the cell membrane, they become able to locally inhibit actin polymerization. Write $c:\R_+ \times \Omega \rightarrow \R$ the concentration in an inactive form, and $\mu: \R_+ \times \partial \Omega \rightarrow \R$ the activated concentration at the boundary. Then, 
\begin{itemize}
\itemr the inactive molecules follow an advection-diffusion dynamics with advection velocity $\bu$,
\itemr there is an exchange dynamics on $\partial \Omega$ between active and inactive forms, 
\itemr $k_p$ is a decreasing function of $\mu$.
\end{itemize}

For $D$ the diffusion coefficient, and $k_{\text{on}/\text{off}}$ the activation/desactivation rates, the corresponding problem writes 

\begin{equation}\label{ReacConvDiff}
\begin{cases}
 \partial_t c(t,\bx) + \Div \left(  c(t,\bx)\bu(t,\bx)  - D \nabla c(t,\bx)\right) = 0 &\, \quad \text{ in } \Omega\, ,\\
 \left(D \nabla c(t,\bx) - c(t,\bx)\bu(t,\bx) \right) \cdot \bn_{\bx} = - k_{\text{on}} c(t,\bx) + k_{\text{off}} \mu(t,\bx) & \, \quad  \text{ on } \partial \Omega\, ,\\
 \frac{\partial}{\partial t} \mu(t,\bx) = k_{\text{on}} c(t,\bx) - k_{\text{off}} \mu(t,\bx)& \, \quad  \text{ on } \partial \Omega\,.
\end{cases}
\end{equation}
Note that the boundary condition ensures mass conservation:
\begin{equation}
\frac{\dx }{\dx t} \left( \int_{\Omega} c(t,\bx) \dx \bx + \int_{\partial \Omega} \mu(t,\bx) \dx \bx \right) =0\,.
\end{equation}

Some remarks can be made to highlight the difficulties in the analysis of the model. The fluid velocity rewrites 
\begin{equation}
\bu(t,\bx) = - \nabla p(t,\bx) - \gamma \int_{\partial \Omega} k_p(t,\bx) \bn_{\bx}\,\dx \bx\,.
\end{equation} 
We notice that this expression depends on the concentration in activated molecules $\mu$ from the pressure boundary term $k_p$, so that the reaction-advection-diffusion problem is nonlinear. Moreover, the integral term makes it also non-local.

The corresponding 1D model in a special case without activation at the boundary and for a linear $k_p$ writes as a nonlinear non-local advection-diffusion problem, and has been studied in \cite{etchegaray2017analysis}. The analysis shows the existence of different asymptotic solutions, describing both motile and non motile behaviours. Moreover, for a subcritical mass of molecules, the global weak existence of solutions is established, along with their convergence to a non motile gaussian profile at explicit rate. Finally, under conditions on the initial condition, it is shown that solutions blow up in finite time. 

Since the model carries non trivial behaviours, computing numerical solutions is of particular interest. In the following, we develop a finite volume method for the polar formulation of the problem \eqref{Poisson}-\eqref{Darcy}-\eqref{ReacConvDiff}, since the domain $\Omega$ is a disc.

\section{Numerical method}
We introduce now the finite volume discretization of the model. We assume that the cytokeleton domain is an annulus $\Omega = B(0,R)\setminus B(0,R_{\min}) \subset \R^2$, where $B(0,R_{\min})$ accounts for the nucleus. We consider a zero-flux boundary condition for the molecular concentration on $C(0,R_{\min})$, the circle of center $(0,0)$ and radius $R_{\min}$. As a consequence, it will be natural to study the problem in polar coordinates.

\subsection{Polar formulation}
\subsubsection*{Reaction-advection-diffusion problem}
 Let $\bx=(r \cos(\theta), r \sin(\theta)) \in \Omega$, and $\tilde{c}$ (similarly $\tilde{\mu}$) the polar function such that $\frac{1}{r} \tilde{c}(t,r,\theta) = c(t,\bx)$ with $(r,\theta) \in [R_{\min},R] \times \mathbb{R}/2\pi\mathbb{Z}$. Then, the problem on the molecular specie writes 
\begin{eqnarray}
\partial_t \tilde{c}(t,r,\theta) & = & \partial_r \left(D r \partial_r \left(\frac{\tilde{c}(t,r,\theta)}{r}\right) - \tilde{c}(t,r,\theta) \bu_r (t,r,\theta) \right) \nonumber \\
& + &  \partial_{\theta} \left( \frac{1}{r^2} (D\partial_{\theta}\tilde{c}(t,r,\theta) - \tilde{c}(t,r,\theta) \bu_\theta (t,r,\theta) ) \right) \label{c1polar} \mbox{ in } \Omega\, , \\
k_{\textnormal{off}} \tilde{\mu}(t,R,\theta)-k_{\textnormal{on}} \tilde{c}(t,R,\theta)  & = &  D R \partial_r \left(\frac{\tilde{c}(t,R,\theta)}{R}\right) - \tilde{c}(t,R,\theta) \bu_r (t,R,\theta) \label{c2polar}, \mbox{ on } C(0,R)\,,  \\
0 &=& D R_{\min} \partial_r \left(\frac{\tilde{c}(t,R_{\min},\theta)}{R_{\min}}\right) - \tilde{c}(t,R_{\min},\theta) \bu_r (t,R_{\min},\theta) ,\nonumber  \mbox{ on } C(0,R_{\min})\,, \\
&&\label{c3polar}\\
\partial_t \tilde{\mu}(t,R,\theta) &=& k_{\textnormal{on}} \tilde{c}(t,R,\theta) -k_{\textnormal{off}} \tilde{\mu}(t,R,\theta)\, , \quad \mbox{ on } C(0,R)\,.
\end{eqnarray}

\subsubsection*{Poisson problem on \texorpdfstring{$p$}{p}}
For the polymerization function $k_p$, let us choose a simple form, that is 
\begin{equation}
k_p(t,\bx) = [1-\delta \mu(t,\bx)]_+\,,
\end{equation}
with $\delta >0$ and $\bx \in \partial \Omega$.\par 
The pressure $p$ is solution of the following problem: 
\begin{eqnarray}
& -\Delta p(t,\bx) = -k_d \label{p1rec}, \mbox{ in } \Omega\ ,\\
& p(t,\bx) = [1 - \delta \mu(t,\bx)]_+ \label{p2rec}, \mbox{ on } C(0,R)\,, \\
&  p(t,\bx) = 0 \label{p3rec}, \mbox{ on } C(0,R_{\min})\,,
\end{eqnarray}
where the pressure condition on $C(0,R_{\min})$ is arbitrary, and fix the pressure values. Let us consider these equations in polar coordinates with $\frac{1}{r} \tilde{p}(t,r,\theta) = p(t,\bx)$ for $(r,\theta) \in [R_{\min},R] \times \mathbb{R}/2\pi\mathbb{Z}$. We have 
\begin{eqnarray}
 - \partial_r \left(r \partial_r \left(\frac{\tilde{p}(r,\theta)}{r}\right) \right)- \frac{1}{r^2} \partial_{\theta \theta} \tilde{p} (r,\theta)  = - k_d \, r, & \quad \mbox{ in } \Omega \label{p1polar},\\
 \tilde{p}(t,R,\theta) =  R\left[1 - \delta \frac{\tilde{\mu}(t,R,\theta)}{R}\right]_+
& \quad  \mbox{ on } C(0,R) \label{p2polar}, \\
 \tilde{p}(t,R_{\min},\theta)= 0 & \quad \mbox{ on } C(0,R_{\min}) \label{p3polar}.
\end{eqnarray}

\subsection{Discretization}
Let $t^n = n\, \Delta t$ be the time discretization, and $\{r_j=R_{\min} + (j-\frac{1}{2}) \, \Delta r, j \in \{1,...,N_r\}\}$ the space discretization of the bounded interval $[R_{\min},R]$, such that $r_{N_r + \frac{1}{2}} = R,$ (therefore $N_R = \frac{R-R_{\min}}{\Delta r}$). Similarly, $\{\theta_k=k \, \Delta \theta, k \in \{1,...,N_\theta\}\}$ is the space discretization of the periodic interval $\mathbb{R}/ 2 \pi \mathbb{Z}$.
 We introduce the control volumes $W_{(j,k)}  \subset \mathbb{R}^2$ and $V_k \subset \R/ 2 \pi \mathbb{Z}$ with
\begin{eqnarray*}
V_k &=& \left(\theta_{k-\frac{1}{2}},\theta_{k+\frac{1}{2}}\right)\,, \\
W_{(j,k)} &=& \left(r_{j-\frac{1}{2}},r_{j+\frac{1}{2}}\right) \times V_k\,.
\end{eqnarray*}

Let $\tilde{c}^{n}_{(j,k)}$ (resp. $\tilde{\mu}^n_{k}$) be the approximated value of the exact solution $\tilde{c}(t^n,r_j,\theta_k)$ (resp. $\tilde{\mu}(t^n,\theta_k)$), and $\tilde{p}^n_{(j,k)}$ be the approximated value of the exact solution $\tilde{p}(t^n,r_j,\theta_k)$:
\begin{eqnarray*}
\tilde{c}^{n}_{(j,k)} &\simeq & \frac{1}{\Delta r \Delta \theta} \iint_{W(j,k)} \tilde{c}(t^n,r,\theta) \dx r \dx \theta \,,\\
\tilde{\mu}^n_{k} &\simeq & \frac{1}{\Delta \theta} \int_{V_k} \tilde{\mu}(t^n,\theta) \dx \theta \,,\\
\tilde{p}^n_{(j,k)} &\simeq & \frac{1}{\Delta r \Delta \theta} \iint_{W(j,k)} \tilde{p}(t^n,r,\theta) \dx r \dx \theta \,.\\
\end{eqnarray*}
Moreover, we write $\bu^n$ and $\bv^n$ the corresponding discretized velocity functions at time $t^n$.

The resolution is made as follows: for $n\geq 0$, knowing $(\tilde{c}^{n},\tilde{\mu}^{n})$ allows to compute $\tilde{p}^n$, then $(\bu^n,\bv^n)$. Finally, $(\tilde{c}^{n+1},\tilde{\mu}^{n+1}) $ is computed using $\bu^n$, and so on. \par 

\subsubsection*{Problem on \texorpdfstring{$\tilde{p}$}{p}}
The pressure problem is time-dependent because of the Dirichlet boundary condition. Therefore, it is solved explicitly in time. Write $\mathcal{F}$ for the numerical flux. Then, we have the following scheme for equation \eqref{p1polar}: for $(j,k) \in \{1,...,N_r\} \times \{1,...,N_\theta\}$,
\begin{eqnarray*}
- \left(\frac{\mathcal{F}_{(j+\frac{1}{2},k)} - \mathcal{F}_{(j-\frac{1}{2},k)}}{\Delta r} +  \frac{\mathcal{F}_{(j,k+\frac{1}{2})} - \mathcal{F}_{(j,k-\frac{1}{2})}}{\Delta \theta}\right)  = -k_d r_j\,.
\end{eqnarray*}
The finite volume numerical fluxes are defined by
\begin{eqnarray*}
\mathcal{F}_{(j+\frac{1}{2},k)} = r_{j+\frac{1}{2}}\, \frac{\frac{\tilde{p}_{(j+1,k)}}{r_{j+1}} - \frac{\tilde{p}_{(j,k)}}{r_{j}} }{\Delta r}, & \quad &
\mathcal{F}_{(j-\frac{1}{2},k)} = r_{j-\frac{1}{2}}\, \frac{\frac{\tilde{p}_{(j,k)}}{r_{j}} - \frac{\tilde{p}_{(j-1,k)}}{r_{j-1}} }{\Delta r},\\
\mathcal{F}_{(j,k+\frac{1}{2})} = \frac{1}{r_{j}^2} \frac{\tilde{p}_{(j,k+1)} - \tilde{p}_{(j,k)} }{\Delta \theta}, & \quad &
\mathcal{F}_{(j,k-\frac{1}{2})} = \frac{1}{r_{j}^2} \frac{\tilde{p}_{(j,k)} - \tilde{p}_{(j,k-1)} }{\Delta \theta}.
\end{eqnarray*}

The Dirichlet boundary conditions \eqref{p2polar}-\eqref{p3polar} are imposed using ghost values $\tilde{p}_{(0,k)}$ and $\tilde{p}^n_{(N_r+1,k)}$. For $k \in \{1,...,N_\theta\}$, 

$$\mathcal{F}_{(\frac{1}{2},k)} = \frac{ r_{\frac{1}{2}}}{r_{1}}\, \frac{\tilde{p}_{(1,k)}}{ \Delta r}\,,$$ 

since $\tilde{p}_{(0,k)}=0$. Now, $\tilde{p}^n_{(N_r+1,k)}= r_{N_r}\left[1 - \delta \frac{\tilde{\mu}^n_{(N_r,k)}}{r_{N_r}}\right]_+$, the corresponding flux writes 
$$\mathcal{F}_{(N_r+\frac{1}{2},k)} = r_{N_r+\frac{1}{2}}\, \frac{ \frac{r_{N_r}}{r_{N_r +1}} \left(1 - \delta \frac{\tilde{\mu}^n_{(N_r,k)}}{r_{N_r}}\right)    - \frac{\tilde{p}_{(N_r,k)}}{r_{N_r}} }{\Delta r}\,.$$

Therefore, the term $ \frac{ r_{N_r} r_{N_r+\frac{1}{2}} }{r_{N_r +1}\Delta r}\left[1 - \delta \frac{\tilde{\mu}^n_{(N_r,k)}}{r_{N_r}}\right]_+ $ will be included in the right hand side of the matricial problem.

Similarly, the periodic conditions impose for $j \in \{1,...,N_r\}$,
$$\mathcal{F}_{(j,N_\theta+\frac{1}{2})} = \mathcal{F}_{(j,\frac{1}{2})}=\frac{1}{r_{j}^2} \frac{\tilde{p}_{(j,1)} - \tilde{p}_{(j,N_\theta)} }{\Delta \theta}.$$

As a consequence, we write the corresponding system, where the terms in bold account for the boundary conditions. Note also that each equation is written for $k\in \{1,...,N_{\theta}\}$ with the convention that if $k=1$, then $k-1=N_{\theta}$ and $k-\frac{1}{2}=\frac{1}{2}$. Similarly, if $k=N_{\theta}$, then $k+1=1$ and $k+\frac{1}{2}=\frac{1}{2}$. 

For $j=1$, we write 

\begin{equation}
\begin{aligned}
&\frac{1}{\Delta r^2} \left[ r_{\frac{3}{2}} \left(\frac{\tilde{p}_{(2,k)}}{r_2} - \frac{\tilde{p}_{(1,k)}}{r_1} \right) - \mathbf{\frac{r_{\frac{1}{2}}}{r_1} \tilde{p}_{(1,k)}} \right] + 
\frac{\tilde{p}_{(1,k+1)} - 2\tilde{p}_{(1,k)} + \tilde{p}_{(1,k-1)} }{r_1^2 \Delta \theta^2}  = k_d r_1\,.\\
\end{aligned}
\end{equation}
or equivalently 

\begin{equation}
\begin{aligned}
&\frac{1}{\Delta r^2} \left[  \frac{r_{\frac{1}{2}}+r_{\frac{3}{2}}}{r_1}\tilde{p}_{(1,k)} - \frac{r_{\frac{3}{2}}}{r_2}\tilde{p}_{(2,k)} \right] + 
\frac{-\tilde{p}_{(1,k-1)}  + 2\tilde{p}_{(1,k)} - \tilde{p}_{(1,k+1)}}{r_1^2 \Delta \theta^2} = - k_d r_1\,.\\
\end{aligned}
\end{equation}

For $j \in \{2,...,N_{\theta}-1\}$, we write for $k\in \{1,...,N_{\theta} \}$,
\begin{equation}
\begin{aligned}
&\frac{1}{\Delta r^2} \left[ r_{j+\frac{1}{2}} \left(\frac{\tilde{p}_{(j+1,k)}}{r_{j+1}} - \frac{\tilde{p}_{(j,k)}}{r_j} \right) - r_{j-\frac{1}{2}} \left(\frac{\tilde{p}_{(j,k)}}{r_j} - \frac{\tilde{p}_{(j-1,k)}}{r_{j-1}} \right) \right] + 
\frac{\tilde{p}_{(j,k+1)} - 2\tilde{p}_{(j,k)} + \tilde{p}_{(j,k-1)} }{r_j^2 \Delta \theta^2} = k_d r_j\,.\\
\end{aligned}
\end{equation}

or equivalently 

\begin{equation}
\begin{aligned}
&\frac{1}{\Delta r^2} \left[ -\frac{r_{j-\frac{1}{2}}}{r_{j-1}} \tilde{p}_{(j-1,k)}
 + \frac{r_{j+\frac{1}{2}} + r_{j-\frac{1}{2}}}{r_j} \tilde{p}_{(j,k)} - \frac{r_{j+\frac{1}{2}}}{r_{j+1}} \tilde{p}_{(j+1,k)} \right] + 
\frac{-\tilde{p}_{(j,k-1)} +2\tilde{p}_{(j,k)}
-\tilde{p}_{(j,k+1)} }{r_j^2 \Delta \theta^2} = -k_d r_j\,.\\
\end{aligned}
\end{equation}

Finally, for $j=N_r$, we write for $k\in \{1,...,N_{\theta} \}$,
\begin{equation}
\begin{aligned}
&
\begin{split}
\frac{1}{\Delta r^2} \left[ \mathbf{-\frac{r_{N_r+\frac{1}{2}}}{r_{N_r}} 
\tilde{p}_{(N_r,k)}} - r_{N_r-\frac{1}{2}} \left(\frac{\tilde{p}_{(N_r,k)}}{r_{N_r}} - \frac{\tilde{p}_{(N_r - 1,k)}}{r_{N_r - 1}} \right) \right] + 
\frac{\tilde{p}_{(N_r,k+1)} - 2\tilde{p}_{(N_r,k)} + \tilde{p}_{(N_r,k-1)} }{r_j^2 \Delta \theta^2} 
 = k_d r_{N_r} \\- \frac{ r_{N_r} r_{N_r+\frac{1}{2}} }{r_{N_r +1}\Delta r} \left[1 - \delta \frac{\tilde{\mu}^n_{N_r,k}}{r_{N_r}}\right]_+ \,.\end{split}\\
\end{aligned}
\end{equation}

or equivalently 

\begin{equation}
\begin{aligned}
&
\begin{split}
\frac{1}{\Delta r^2} \left[-\frac{r_{N_r-\frac{1}{2}}}{r_{N_r - 1}} \tilde{p}_{(N_r - 1,k)}
+\frac{r_{N_r-\frac{1}{2}}+r_{N_r+\frac{1}{2}}}{r_{N_r}} 
\tilde{p}_{(N_r,k)} \right] + 
\frac{-\tilde{p}_{(N_r,k-1)} + 2\tilde{p}_{(N_r,k)} - \tilde{p}_{(N_r,k+1)} }{r_j^2 \Delta \theta^2}  = -k_d r_{N_r} \\
+ \frac{ r_{N_r} r_{N_r+\frac{1}{2}} }{r_{N_r +1}\Delta r} \left[1 - \delta \frac{\tilde{\mu}^n_{N_r,k}}{r_{N_r}}\right]\,,
\end{split}
\\
\end{aligned}
\end{equation}

Let us now write the corresponding matricial problem. We define the column vector $\mathcal{P}$ by $\mathcal{P}(k+(j-1) N_\theta) = \tilde{p}_{(j,k)}$ with $(j,k) \in \{1,...,N_r\} \times \{1,...,N_\theta\}$:
\begin{equation*}
\mathcal{P} =
\left(
\tilde{p}_{(1,1)} \, \dots \,  \tilde{p}_{(1,N_\theta)} \, \tilde{p}_{(2,1)} \, \dots \, \tilde{p}_{(2,N_\theta)} \, \dots  \, \tilde{p}_{(N_r,N_\theta)}
\right)^T
\end{equation*}
For $\Delta r=\Delta \theta$ the stiffness matrix $\mathcal{A}_p$ is defined by
\begin{eqnarray}\label{A2D}
\mathcal{A}_p = 
\begin{pmatrix} 
&\frac{r_{1/2}+r_{1+ 1/2} }{r_1} I_{N_\theta} & -\frac{r_{1+\frac{1}{2}}}{r_{2}} \, I_{N_\theta} \\
\\
& \ddots & \ddots & \ddots &  \\ 
& & -\frac{r_{j-\frac{1}{2}}}{r_{j-1}} \,I_{N_\theta} & \frac{r_{j-\frac{1}{2}} + r_{j+\frac{1}{2}}}{r_j} \,I_{N_\theta} & -\frac{r_{j+\frac{1}{2}}}{r_{j+1}} \, I_{N_\theta} \\ & & & \ddots & \ddots & \ddots &  \\ 
\\
&&&& -\frac{r_{N_r-\frac{1}{2}}}{r_{N_r-1}} \,I_{N_\theta}  & \frac{r_{N_r-\frac{1}{2}}+r_{N_r+\frac{1}{2}}}{r_{N_r}} I_{N_\theta} 
\end{pmatrix} \nonumber \\
+ 
\begin{pmatrix} 
\frac{1}{r_{1}^2} A  \\  & \frac{1}{r_{2}^2} A  \\  & & \ddots \\ & & & \frac{1}{r_{N_r-1}^2} A \\ & & & & \frac{1}{r_{N_r}^2} A
\end{pmatrix},
\end{eqnarray}
where the second matrix accounts for the angular diffusion, with $A \in M_{N_\theta} (\mathbb{R})$ the classical diffusion matrix with periodic flux boundary conditions:
\begin{equation}\label{defApolar}
A = \begin{pmatrix} 
2 & -1 & & & -1 \\ -1 & 2 & \ddots \\ & \ddots & \ddots & \ddots \\& & \ddots & 2 & -1 \\ -1 & & & -1 & 2
\end{pmatrix}\,.
\end{equation}

The right-hand side in \eqref{p1polar} and the flux boundary condition \eqref{p2polar} on $C(0,R_{\max})$ imposes this right hand side column vector of length $N_r \, N_\theta$:
\begin{equation*}
\mathcal{R}^n_p = -k_d \begin{pmatrix}
\left(r_{1}\right)_{k} \\ \vdots \\ \left(r_{j}\right)_{k} \\ \vdots \\ \left(r_{N_r}\right)_{k} 
\end{pmatrix}
+ \frac{ r_{N_r} r_{N_r+\frac{1}{2}} }{r_{N_r +1}\Delta r} 
 \begin{pmatrix} 0 \\ \vdots \\ 0 \\ {\left[1 - \delta \frac{\tilde{\mu}^n_{N_r,k}}{r_{N_r}}\right]_+} \end{pmatrix} \,.
\end{equation*}

We use a standard numerical method to invert the symmetric positive definite matrix $\frac{1}{\Delta r^2} \mathcal{A}_p $ and then resolve at each time step
\begin{equation*}
\mathcal{P} = \left(\frac{1}{\Delta r^2}\mathcal{A}_p
\right)^{-1} \,  \mathcal{R}^n_p.
\end{equation*}

\subsubsection*{Equations for \texorpdfstring{$\bu$ and $\bv$}{u and v}}
The velocities $\bu$ and $\bv$ depend on the pressure $p$, which is obtained through a time explicit scheme. Therefore, they also depend explicitly in time on the concentration.\par 
The equation on $\bv$ in polar coordinates writes 
\begin{equation}
\bv(t) =  \gamma \int_{0}^{2\pi} \left[1-\delta \frac{\tilde{\mu}(t,R,\theta)}{R}\right]_+ \bn \dx\theta \,.
\end{equation}
We compute numerically the velocity in cartesian coordinates $\bv^n_{\textnormal{cart}} := (v_x^n,v_y^n)^T$:
\begin{eqnarray}
v_x^n &=&  \gamma \Delta \theta \sum_{k=1}^{N_\theta} \left[1- \delta \frac{\tilde{\mu}^{n}_{(N_r,k)}}{R}\right]_+ \cos(\theta_k)\,,\\
v_y^n &=&  \gamma \Delta \theta \sum_{k=1}^{N_\theta} \left[1- \delta \frac{\tilde{\mu}^{n}_{(N_r,k)}}{R}\right]_+ \sin(\theta_k)\,.
\end{eqnarray}
Then, a polar change of coordinates leads to
$\bv^n := (v_r^n,v_{\theta}^n)^T$.

For the fluid velocity, we have
\begin{equation}
\bu(t,\bx) = - \nabla p(t,\bx)  - \bv(t) \,,
\end{equation}
that rewrites
\begin{equation}\label{eq:u_rigide_pol}
\bu(t,r,\theta) = - \left( \partial_r \left(\frac{\tilde{p}(t,r,\theta)}{r} \right) +v_r \right) \vec{\be_r}
- \left(  \frac{1}{r} \partial_\theta \tilde{p}(t,r,\theta) +  v_{\theta}\right) \vec{\be_{\theta}}  \,,
\end{equation}
since $u_r = u(t,\bx) \cdot \vec{\be_r}$ and $u_\theta = r u(t,\bx) \cdot \vec{\be_\theta}$. We define at time $t^n$
\begin{eqnarray*}
u^n_{(j+\frac{1}{2},k)}= - \frac{\frac{\tilde{p}_{(j+1,k)}}{r_{j+1}} - \frac{\tilde{p}_{(j,k)}}{r_{j}} }{\Delta r} - v_r^n\,, & \quad &
u^n_{(j-\frac{1}{2},k)}= - \frac{\frac{\tilde{p}_{(j,k)}}{r_{j}} - \frac{\tilde{p}_{(j-1,k)}}{r_{j-1}}}{\Delta r} - v_r^n\,, \\
u^n_{(j,k+\frac{1}{2})}=- \frac{1}{r_j}\frac{\tilde{p}_{(j,k+1)} - \tilde{p}_{(j,k)} }{\Delta \theta} -v_\theta^n\,, & \quad &
u^n_{(j,k-\frac{1}{2})}=- \frac{1}{r_j}\frac{\tilde{p}_{(j,k)} - \tilde{p}_{(j,k-1)} }{\Delta \theta}-v_\theta^n\,.
\end{eqnarray*}

\subsubsection{Problem for \texorpdfstring{$\tilde{c}$}{c} and \texorpdfstring{$\tilde{\mu}$}{u}}
For simplicity, we call  again $\mathcal{F}$ the numerical fluxes. We can write the following scheme for equation \eqref{c1polar}: for $(j,k) \in \{1,...,N_r\} \times \{1,...,N_\theta\}$, 
\begin{eqnarray*}
\frac{\tilde{c}_{(j,k)}^{n+1} - \tilde{c}_{(j,k)}^n}{\Delta t} = \frac{\mathcal{F}_{(j+\frac{1}{2},k)} - \mathcal{F}_{(j-\frac{1}{2},k)}}{\Delta r} +  \frac{\mathcal{F}_{(j,k+\frac{1}{2})} - \mathcal{F}_{(j,k-\frac{1}{2})}}{\Delta \theta}.
\end{eqnarray*}
We define now the numerical fluxes. The diffusion part is implicit, so that no CFL condition is needed \citep{allaire2005analyse}, while the advection is explicit due to the nonlinearity in the expression of $v$.

We have:
\begin{eqnarray*}
\mathcal{F}_{(j+\frac{1}{2},k)} &=& D\, r_{j+\frac{1}{2}}\, \frac{\frac{\tilde{c}^{n+1}_{(j+1,k)}}{r_{j+1}} - \frac{\tilde{c}^{n+1}_{(j,k)}}{r_{j}} }{\Delta r}  - A^{up} \left(u^{n}_{(j+\frac{1}{2},k)},\tilde{c}^{n}_{(j,k)},\tilde{c}^{n}_{(j+1,k)}\right),\\
\mathcal{F}_{(j-\frac{1}{2},k)} &=& D\, r_{j-\frac{1}{2}}\, \frac{\frac{\tilde{c}^{n+1}_{(j,k)}}{r_{j}} - \frac{\tilde{c}^{n+1}_{(j-1,k)}}{r_{j-1}} }{\Delta r}  - A^{up} \left(u^{n}_{(j-\frac{1}{2},k)},\tilde{c}^{n}_{(j-1,k)},\tilde{c}^{n}_{(j,k)}\right),\\
\mathcal{F}_{(j,k+\frac{1}{2})}  &=& \frac{1}{r_{j}^2} \, \left(D\, \frac{\tilde{c}^{n+1}_{(j,k+1)} - \tilde{c}^{n+1}_{(j,k)} }{\Delta \theta} - A^{up} \left(u^{n}_{(j,k+\frac{1}{2})},\tilde{c}^{n}_{(j,k)},\tilde{c}^{n}_{(j,k+1)}\right) \right),\\
\mathcal{F}_{(j,k+\frac{1}{2})} &=&  \frac{1}{r_{j}^2} \, \left(D\, \frac{\tilde{c}^{n+1}_{(j,k)} - \tilde{c}^{n+1}_{(j,k-1)} }{\Delta \theta} - A^{up} \left(u^{n}_{(j,k-\frac{1}{2})},\tilde{c}^{n}_{(j,k-1)},\tilde{c}^{n}_{(j,k)}\right)\right)\,,
\end{eqnarray*}
where $\textit{A}^{up}$ is the advection term expressed by

\begin{equation}\label{Auppolar}
\textit{A}^{up}(\bu,x_-,x_+) = \begin{cases}
\bu \, x_-, \quad \mbox{ if } u > 0, \\ 
\bu \, x_+, \quad \mbox{ if } u < 0.
\end{cases}
\end{equation}
The external boundary condition \eqref{c2polar} yields 
$$\mathcal{F}_{(N_r + \frac{1}{2},k)} = k_{\textnormal{off}}\, \mu^{n+1}_{k} - k_{\textnormal{on}}\, \tilde{c}^{n+1}_{N_r,k}$$ for $k \in \{1,...,N_\theta\}$.
The zero flux boundary condition \eqref{c3polar} imposes that $\mathcal{F}_{(\frac{1}{2},k)} = 0$ for $k \in \{1,...,N_\theta\}$. 
Similarly, the periodic conditions impose for $j \in \{1,...,N_r\}$ 
$$\mathcal{F}_{(j,N_\theta+\frac{1}{2})} = \mathcal{F}_{(j,\frac{1}{2})}=\frac{1}{r_{j}^2} \, \left( D\, \frac{\tilde{c}^{n+1}_{(j,1)} - \tilde{c}^{n+1}_{(j,N_\theta)} }{\Delta \theta} - A^{up} \left(u^{n}_{(j,\frac{1}{2})},\tilde{c}^{n}_{(j,N_\theta)},\tilde{c}^{n}_{(j,1)}\right)\right).$$

We write the corresponding scheme and group the implicit (resp. explicit) terms on the left-hand-side (resp. right-hand-side). Note also that each equation is written for $k\in \{1,...,N_{\theta}\}$ with the convention that if $k=1$, then $k-1=N_{\theta}$ and $k-\frac{1}{2}=\frac{1}{2}$. Similarly, if $k=N_{\theta}$, then $k+1=1$ and $k+\frac{1}{2}=\frac{1}{2}$. 

For $j=1$, we have
\begin{equation}
\begin{aligned}
&\begin{split}
\left(1+\frac{D \Delta t}{\Delta r^2}  \frac{r_{1+\frac{1}{2}}}{r_1}  \right)\co{n+1}{1}{k}
- \frac{D \Delta t}{\Delta r^2} \frac{r_{1+\frac{1}{2}}}{r_{2}} \co{n+1}{2}{k} 
+ \frac{D \Delta t}{r_1^2 \Delta \theta^2} \left( - \co{n+1}{1}{k-1} +2 \co{n+1}{1}{k} - \co{n+1}{1}{k+1} \right) \\
 = \co{n}{1}{k} - \frac{\Delta t}{ \Delta r} A^{up} \left(u^{n}_{(1+\frac{1}{2},k)},\tilde{c}^{n}_{(1,k)},\tilde{c}^{n}_{(2,k)}\right)\\
+ \frac{\Delta t}{r_1^2 \Delta \theta} \left( A^{up} \left(u^{n}_{(1,k-\frac{1}{2})},\tilde{c}^{n}_{(1,k-1)},\tilde{c}^{n}_{(1,k)}\right) - A^{up} \left(u^{n}_{(1,k+\frac{1}{2})},\tilde{c}^{n}_{(1,k)},\tilde{c}^{n}_{(1,k+1)}\right) \right).\end{split}
\\
\end{aligned}
\end{equation}

Now, for $j\in{2,...,N_r-1}$, we obtain
\begin{equation}
\begin{aligned}
&\begin{split}
\co{n+1}{j}{k} + \frac{D \Delta t}{\Delta r^2} \left(- \frac{r_{j-\frac{1}{2}}}{r_{j-1}} \co{n+1}{j-1}{k} + \frac{r_{j-\frac{1}{2}} + r_{j+\frac{1}{2}}}{r_j} \co{n+1}{j}{k} - \frac{r_{j+\frac{1}{2}}}{r_{j+1}} \co{n+1}{j+1}{k} \right)
+ \frac{D \Delta t}{r_j^2 \Delta \theta^2} \left( - \co{n+1}{j}{k-1} +2 \co{n+1}{j}{k} - \co{n+1}{j}{k+1} \right) \\
 = \co{n}{j}{k} + \frac{\Delta t}{ \Delta r} \left(  A^{up} \left(u^{n}_{(j-\frac{1}{2},k)},\tilde{c}^{n}_{(j-1,k)},\tilde{c}^{n}_{(j,k)}\right) - A^{up} \left(u^{n}_{(j+\frac{1}{2},k)},\tilde{c}^{n}_{(j,k)},\tilde{c}^{n}_{(j+1,k)}\right)\right)\\
+ \frac{\Delta t}{r_j^2 \Delta \theta} \left( A^{up} \left(u^{n}_{(j,k-\frac{1}{2})},\tilde{c}^{n}_{(j,k-1)},\tilde{c}^{n}_{(j,k)}\right) - A^{up} \left(u^{n}_{(j,k+\frac{1}{2})},\tilde{c}^{n}_{(j,k)},\tilde{c}^{n}_{(j,k+1)}\right) \right).\end{split}
\\
\end{aligned}
\end{equation}

Finally, for $j=N_r$, we get
\begin{equation}
\begin{aligned}
&\begin{split}
- \frac{k_{\text{off}}\Delta t}{\Delta r} \mu_k^{n+1}  \left(1+ \frac{k_{\text{on}}\Delta t}{\Delta r} 
+ \frac{D \Delta t}{\Delta r^2}\frac{r_{N_r-\frac{1}{2}}}{r_{N_r}} \right)\co{n+1}{N_r}{k} - \frac{D \Delta t}{\Delta r^2} \frac{r_{N_r-\frac{1}{2}}}{r_{N_r-1}} \co{n+1}{N_r-1}{k} 
\\ + \frac{D \Delta t}{r_{N_r}^2 \Delta \theta^2} \left( - \co{n+1}{{N_r}}{k-1} +2 \co{n+1}{{N_r}}{k} - \co{n+1}{{N_r}}{k+1} \right) 
 = \co{n}{{N_r}}{k} + \frac{\Delta t}{ \Delta r} \left(  A^{up} \left(u^{n}_{({N_r}-\frac{1}{2},k)},\tilde{c}^{n}_{({N_r}-1,k)},\tilde{c}^{n}_{({N_r},k)}\right)\right)\\
+ \frac{\Delta t}{r_{N_r}^2 \Delta \theta} \left( A^{up} \left(u^{n}_{({N_r},k-\frac{1}{2})},\tilde{c}^{n}_{({N_r},k-1)},\tilde{c}^{n}_{({N_r},k)}\right) - A^{up} \left(u^{n}_{({N_r},k+\frac{1}{2})},\tilde{c}^{n}_{({N_r},k)},\tilde{c}^{n}_{({N_r},k+1)}\right) \right).\end{split}
\\
\end{aligned}
\end{equation}

The membrane activation equation writes in polar coordinates
\begin{eqnarray}\label{mupolar}
\partial_t \mu(t,R,\theta) = k_{\textnormal{on}} c(t,R,\theta) -k_{\textnormal{off}} \mu(t,R,\theta)\, , \quad \mbox{ on } C(0,R)\,.
\end{eqnarray}
At each time step, the implicit discretization of equation \eqref{mupolar} for $k \in \{1,...,N_{\theta}\}$ writes
\begin{equation}\label{mu_discretise}
- \Delta t k_{\textnormal{on}} \tilde{c}^{n+1}_k +
(1 + \Delta t \,  k_{\textnormal{off}}) \mu^{n+1}_k  = \mu^n_k\,.
\end{equation}

For simplicity, we treat both the free and activated concentrations in the same linear problem of unknown
\begin{displaymath}
\mathcal{E}^n =
\left( 
\tilde{c}^n_{(1,1)} \, \dots  \, \tilde{c}^n_{(1,N_\theta)} \,  \tilde{c}^n_{(2,1)} \, \dots \, \tilde{c}^n_{(2,N_\theta)} \, \dots  \, \tilde{c}^n_{(N_r,N_\theta)} \, \tilde{\mu}_1^n \, \dots \, \tilde{\mu}_{N_{\theta}}^n 
\right)^T.
\end{displaymath}

We have at each time step
\begin{equation*}
\left( I_{(N_r + 1)N_\theta} + \frac{\Delta t}{\Delta r^2} \mathcal{A}\right) \mathcal{E}^{n+1} = \left( I_{(N_r + 1)N_\theta} - \frac{\Delta t}{\Delta r} \mathcal{B}^{n} \  \right)\mathcal{E}^{n} \,,
\end{equation*}
where $\mathcal{A} $ is the diffusion matrix, and $\mathcal{B}^{n} $ the advection matrix.


The diffusion matrix writes 
\begin{eqnarray}\label{A2DC}
\mathcal{A} = 
\begin{pmatrix} 
&\frac{r_{1/2} }{r_1} I_{N_\theta} & -\frac{r_{1+\frac{1}{2}}}{r_{2}} \, I_{N_\theta}& \\
\\
& \ddots & \ddots & \ddots &  &\\ 
& & -\frac{r_{j-\frac{1}{2}}}{r_{j-1}} \,I_{N_\theta} & \frac{r_{j-\frac{1}{2}} + r_{j+\frac{1}{2}}}{r_j} \,I_{N_\theta} & -\frac{r_{j+\frac{1}{2}}}{r_{j+1}} \, I_{N_\theta} \\ & & & \ddots & \ddots & \ddots &  &\\ 
\\
&&&& -\frac{r_{N_r-\frac{1}{2}}}{r_{N_r-1}} \,I_{N_\theta}  & \left( \frac{r_{N_r+\frac{1}{2}}}{r_{N_r}}+ \Delta r k_{\text{on}} \right)I_{N_\theta} & -\Delta r k_{\text{off}} I_{N_\theta} \\
&&&&& -\Delta r^2 k_{\text{on}}I_{N_\theta} & \Delta r^2 k_{\text{off}}I_{N_\theta}
\end{pmatrix} \nonumber \\
+ 
\begin{pmatrix} 
\frac{1}{r_{1}^2} A  \\  & \frac{1}{r_{2}^2} A  \\  & & \ddots \\ & & & \frac{1}{r_{N_r-1}^2} A \\ & & & & \frac{1}{r_{N_r}^2} A \\
&&&&& 0_{N_{\theta}}
\end{pmatrix},
\end{eqnarray}
where the second matrix accounts for the angular diffusion, with $A \in M_{N_\theta} (\mathbb{R})$ the classical diffusion matrix with periodic flux boundary conditions:
\begin{equation}\label{defApolarC}
A = \begin{pmatrix} 
2 & -1 & & & -1 \\ -1 & 2 & \ddots \\ & \ddots & \ddots & \ddots \\& & \ddots & 2 & -1 \\ -1 & & & -1 & 2
\end{pmatrix}\,.
\end{equation}

Now, for the advection term $A^{up}$ defined by equation \eqref{Auppolar}, we write $(u)^+=\max(u,0)$ and $(u)^-=\min(u,0)$ so that $\textit{A}^{up}(u,\tilde{c}_{j,k},\tilde{c}_{j,k+1}) = (u)^+ c_{j,k} + (u)^- c_{j,k+1}$. Therefore, we introduce the following diagonal matrices for $j \in \{1,...,N_r\}$,  $U^{+}_{j+\frac{1}{2}} \in M_{N_\theta}(\mathbb{R})$ and $U^{-}_{j+\frac{1}{2}} \in M_{N_\theta}(\mathbb{R})$:
\begin{eqnarray*}
U^{\pm}_{j+\frac{1}{2}} = 
\begin{pmatrix} 
\ddots & & & &  \\ & (u_{(j+\frac{1}{2},k-1)}^n)^{\pm} &   \\ & & (u_{(j+\frac{1}{2},k)}^n)^{\pm}  &  &  \\ & & & (u_{(j+\frac{1}{2},k+1)}^n)^{\pm} &   \\  & & & & \ddots
\end{pmatrix} \,.\\
\end{eqnarray*}
Thus we can define the advection matrix:
\begin{eqnarray}\label{B2D}
\mathcal{B}^n & =& 
\begin{pmatrix} 
U_{\frac{3}{2}}^+ & U_{\frac{3}{2}}^- & & &   \\ & \ddots & \ddots  \\ & & U_{j+\frac{1}{2}}^+ & U_{j+\frac{1}{2}}^-  \\ & & & \ddots & U_{N_r-\frac{1}{2}}^- \\ & & & & 0\\
&&&&&0_{N_{\theta}}
\end{pmatrix}
-
\begin{pmatrix} 
0 & & & &  \\ U_{\frac{3}{2}}^+ & \ddots &   \\ & U_{j-\frac{1}{2}}^+ & U_{j-\frac{1}{2}}^- &  \\ & & \ddots & \ddots &  \\  & & & U_{N_r-\frac{1}{2}}^+ & U_{N_r-\frac{1}{2}}^-\\
&&&&&0_{N_{\theta}}
\end{pmatrix} \nonumber \\
& + &
\begin{pmatrix} 
\frac{1}{r_1^2} B^n  \\  & \frac{1}{r_2^2} B^n  \\  & & \ddots \\ & & & \frac{1}{r_{N_r-1}^2}  B^n \\ & & & &  \frac{1}{r_{N_r}^2} B^n\\
&&&&&0_{N_{\theta}}
\end{pmatrix}
\end{eqnarray}
where the discrete advection matrix $B^n \in M_{N_\theta} (\mathbb{R})$ with periodic flux condition on the boundary is defined as in \cite{allaire2005analyse}
\begin{equation}\label{defBpolar}
B^n = 
\begin{pmatrix} 
\left(u_{\frac{3}{2}}^n\right)^+ & \left(u_{\frac{3}{2}}^n\right)^- & & &   \\ & \ddots & \ddots  \\ & & \left(u_{j+\frac{1}{2}}^n\right)^+ & \left(u_{j+\frac{1}{2}}^n\right)^-  \\ & & & \ddots & \left(u_{N_\theta-\frac{1}{2}}^n\right)^- \\ \left(u_{N_\theta+\frac{1}{2}}^n\right)^- & & & & \left(u_{N_\theta+\frac{1}{2}}^n\right)^+ 
\end{pmatrix}
\end{equation}
\begin{equation*}
-
\begin{pmatrix} 
\left(u_{\frac{1}{2}}^n\right)^- & & & & \left(u_{\frac{1}{2}}^n\right)^+  \\ \left(u_{\frac{3}{2}}^n\right)^+ & \ddots &   \\ & \left(u_{j-\frac{1}{2}}^n\right)^+ & \left(u_{j-\frac{1}{2}}^n\right)^- &  \\ & & \ddots & \ddots &  \\  & & & \left(u_{N_\theta-\frac{1}{2}}^n\right)^+ & \left(u_{N_\theta-\frac{1}{2}}^n\right)^-
\end{pmatrix}.
\end{equation*}

\section{Results}
The discretization scheme was implemented using MATLAB. We performed some numerical simulations to test the scheme's numerical convergence. Since the problem has a boundary nonlinearity, comparing a numerical solution to an exact one is out of reach in general. 

We consider the case of a polarised initial condition that writes $\tilde{c}(0,r,\theta)=\cos(\theta-\pi)+1$, $\tilde{\mu}(0,\theta)=0.5 \tilde{c}(0,R,\theta)$. We fix some parameters: $k_{\text{off}}=D=k_d=1$ ; $\delta = \gamma = 2$.

\subsection*{Illustrative example}
\begin{figure}[H]
\centering
\subfloat[$k_{\text{on}}=0.3$]{\includegraphics[scale=0.45]{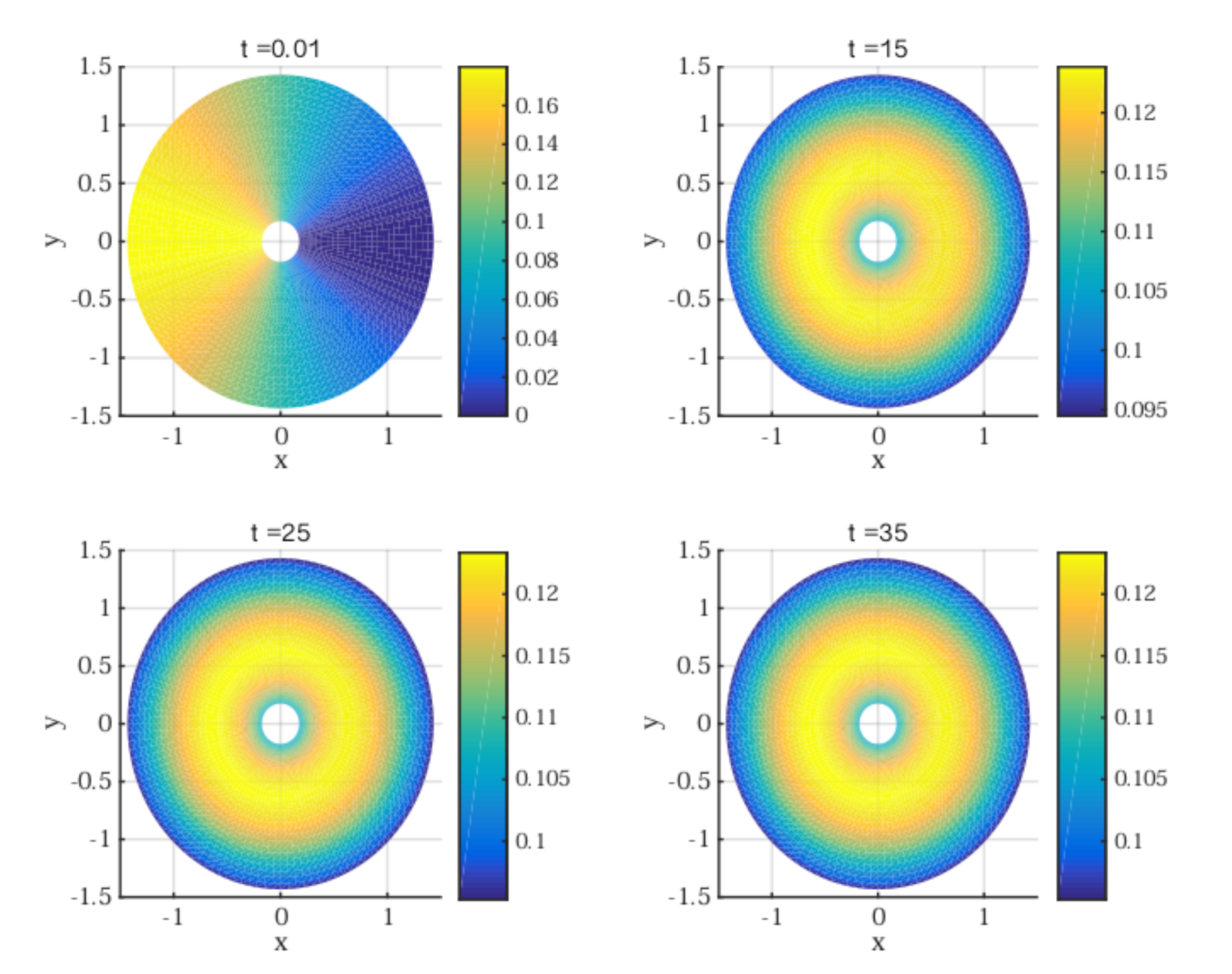}}
\subfloat[$k_{\text{on}}=3$]{\includegraphics[scale=0.45]{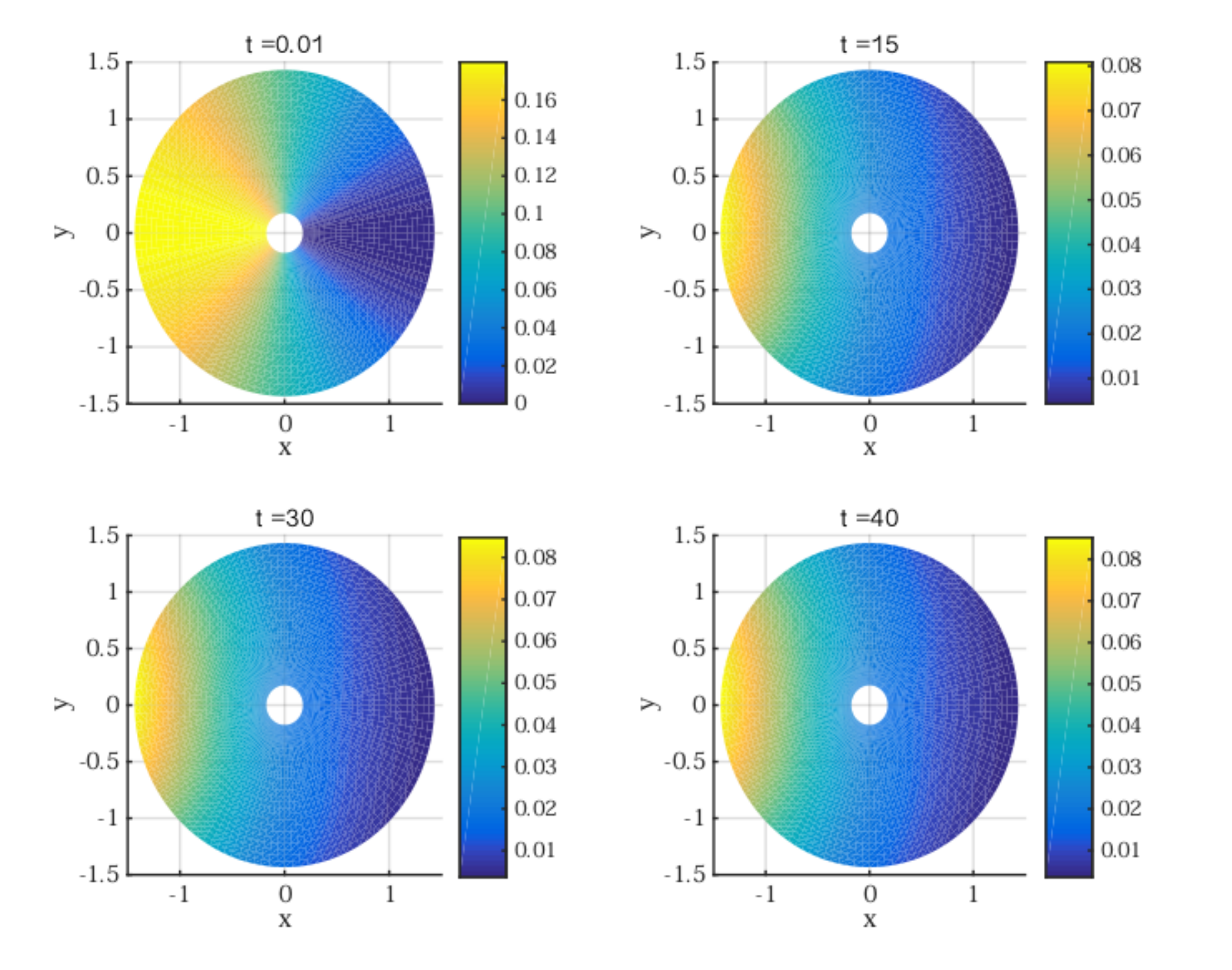}}
\caption{Numerical simulation of the evolution in the molecular concentration over time. Discretization parameters: $\Delta t = 10^{-2}$, $\Delta r = \Delta \theta = 2\pi/120 \simeq 5.3 * 10^{-2}$. Parameters: $R=1.5$, $k_{\text{off}}=D=k_d=1$ ; $\delta = \gamma = 2$. Initial condition: $\tilde{c}(0,r,\theta)=\cos(\theta-\pi)+1$, $\tilde{\mu}(0,\theta)=0.5 \tilde{c}(0,R,\theta)$.}\label{fig:results}
\end{figure}

By changing the value of the activation rate $k_{\text{on}}$, we observe qualitatively different stationary solutions. The figure \ref{fig:results} shows the time evolution in the molecular concentration in the cell body in two characteristic cases. The left figure shows a typically non motile profile, while the right figure displays a polarisation situation.

Note that for $k_{\text{on}}=0$, it is clear that the stationary concentration in activated molecules is $0$ everywhere on the boundary, so that $\bv=0$ and the nonlinearity disappears. In the following, we rather focus on cases where the stationary state can be asymmetric.

\subsection*{Non-zero stationary polarisation}
We also performed the same simulation for $k_{\text{on}}=0.3$ and $R=1$, with an angular discretization step $\Delta \theta = 2\pi / 160 \simeq 3.9 * 10^{-2}$, and varying values for $\Delta r$ and $\Delta t$. The system is considered at numerical steady state when the concentration $\mu$ in activated molecules has stayed unchanged for $1$ numerical hour. Then, we obtain the time to attain the steady state, as well as the stationary polarisation of the system, quantified up to a constant by the cell velocity vector. \par 

We consider the following parameter values:\\
\begin{center}
\begin{tabular}{c|ccccc}
$\Delta r$ & $5*10^{-3}$ &$10^{-2}$ & $2*10^{-2}$ & $2.5 * 10^{-2}$& \\ 
$\Delta t$&$10^{-3}$ & $2*10^{-3}$ & $5*10^{-3}$ & $10^{-2}$ &$1.5 * 10^{-2}$
\end{tabular}
\end{center}

\begin{figure}
\centering
\includegraphics[scale=0.5]{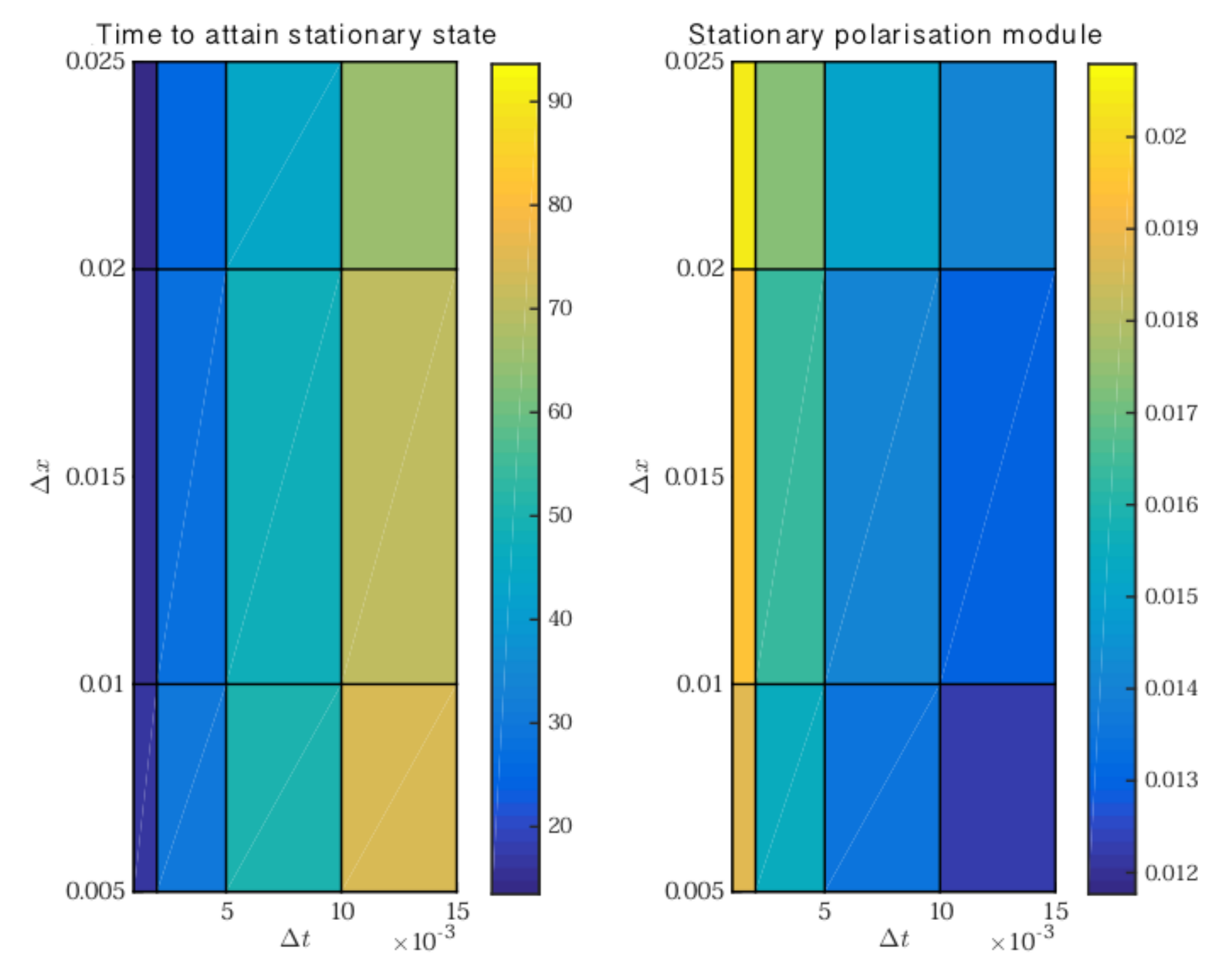}
\caption{Left: Time to attain the stationary state and Right: stationary polarisation module for varying $\Delta r$ and $\Delta t$, and $k_{\text{on}}=0.3$. }\label{fig:surfKon03}
\end{figure}

The figure \ref{fig:surfKon03} shows that the time to attain the stationary state is a consequence of the polarisation state of the stationary solution. Indeed, considering that the initial condition is polarised, and the simulation shows the system's depolarisation towards a less polarised steady state, then the more polarised it is, the sooner it is attained. \par 
We also clearly notice that the smaller $\Delta t$ gets, the more polarised the stationary state gets. This trend can be distinguished independently of $\Delta r$. The radial step has a smaller effect, but more precise grids are correlated with lower polarised solutions. \par 
Overall, the simulations show that a fine time step is fundamental to catch the relevant stationary type of behaviour. The polarisation module being decreasing with the time step, we can infer that the numerical solution approach a polarised state. However, this also shows that the previous numerical simulations are better understood in a qualitative sense rather than for quantitative purposes.

\subsection*{Low stationary polarisation}
\begin{figure}
\centering
\includegraphics[scale=0.5]{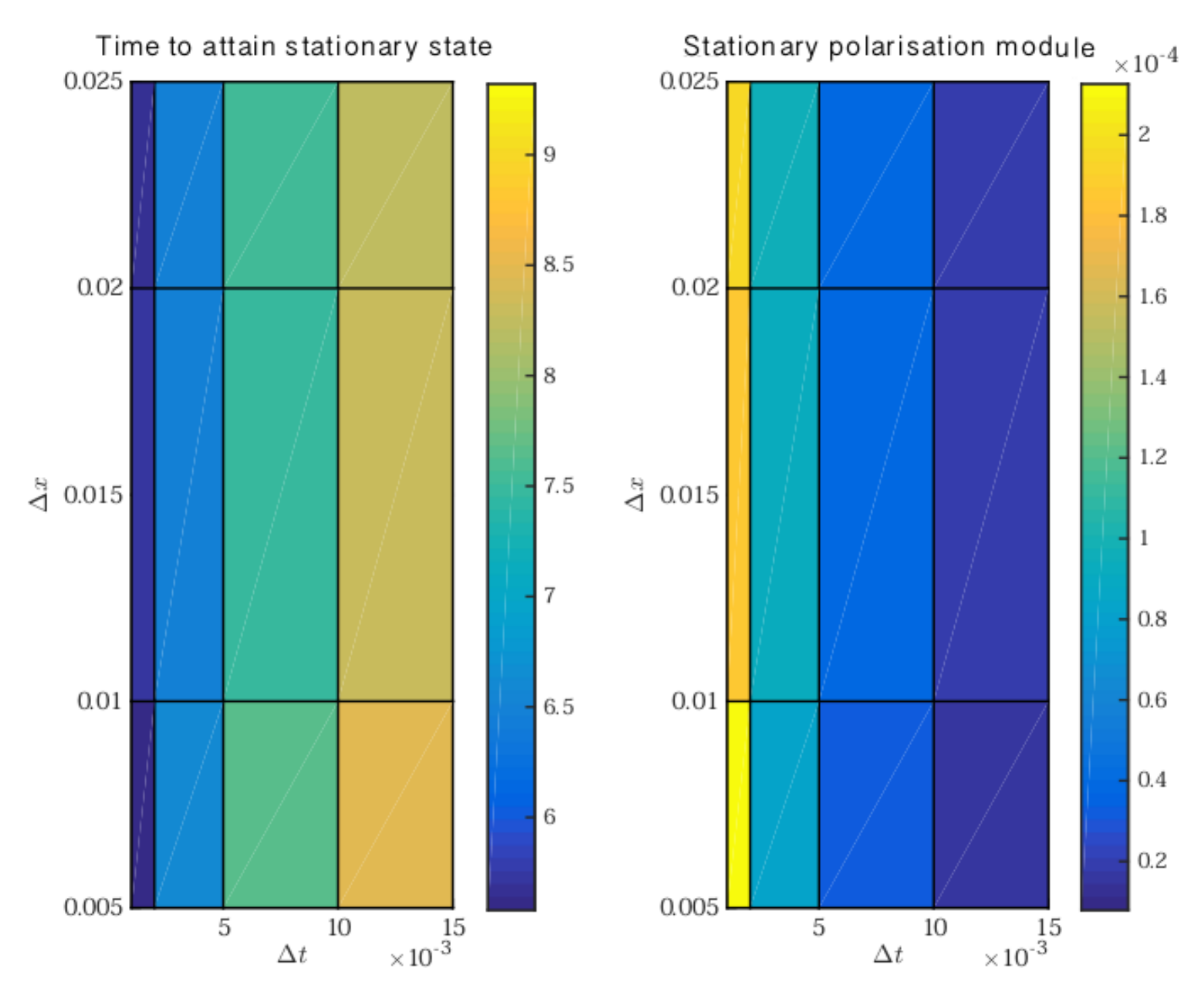}
\caption{Left: Time to attain the stationary state and Right: stationary polarisation module for varying $\Delta r$ and $\Delta t$ and $k_{\text{on}}=0.1$. }\label{fig:surfKon01}
\end{figure}

Finally, we perform the same numerical test with $k_{\text{on}}=0.1$ to check how the discretization steps may generate an error between a polarised and a non polarised stationary solution. We obtain the values depicted in figure \ref{fig:surfKon01}.

We can see that the stationary polarisation levels are very low compared to the case where $k_{\text{on}}=0.3$. The same trend appears for the effect of $\Delta t$, while the effect of $\Delta r$ is less visible since the nonlinear term is small. \par 
To determine if these polarisation levels could be comparable to a true symmetric stationary state, we performed the same simulation for $k_{\text{on}}=0$ and for the most precise time-space grid. The stationary polarisation module was approximately equal to $1.63*10^{-4}$, so that no clear distinction can be made between these cases.

\section{Conclusion}
In this work, we have presented the finite volume discretization of a multiscale model for 2D cell crawling migration consisting in a Darcy fluid dynamics coupled to a Poisson problem and to a nonlinear and non-local reaction-advection-diffusion problem for the concentration in a molecular specie. The simulation of numerical solutions of this type of problems is very useful since the mathematical analysis is necessarily limited, whereas the model show varied behaviours. \par 
The discretization method showed good qualitative numerical result. In particular, the molecular mass is preserved numerically. However, in critical cases, the scheme seems not able to make the distinction between polarised and unpolarised states. 
A natural continuation will consist in taking into account the interaction with the environment (chemical signal, mechanical obstacle) as a bias for motion. Finally, further studies should be based on an implicit treatment of the nonlinearity \cite{cances2017numerical,cances2017nonlinear}.

\bibliographystyle{apalike}
\bibliography{Bib_Asud}

\begin{thebibliography}{}

\bibitem[Allaire, 2005]{allaire2005analyse}
Allaire, G. (2005).
\newblock {\em Analyse num{\'e}rique et optimisation: Une introduction {\`a} la
  mod{\'e}lisation math{\'e}matique et {\`a} la simulation num{\'e}rique}.
\newblock Editions Ecole Polytechnique.

\bibitem[Blanch-Mercader and Casademunt,
  2013]{Blanch-Mercader2013Spontaneous-Mot}
Blanch-Mercader, C. and Casademunt, J. (2013).
\newblock Spontaneous motility of actin lamellar fragments.
\newblock {\em Physical Review Letters}, 110(7).

\bibitem[Bravo-Cordero et~al., 2012]{Bravo-Cordero2012Directed-cell-i}
Bravo-Cordero, J.~J., Hodgson, L., and Condeelis, J. (2012).
\newblock Directed cell invasion and migration during metastasis.
\newblock {\em Curr Opin Cell Biol}, 24(2):277--83.

\bibitem[Canc{\`e}s et~al., 2017a]{cances2017nonlinear}
Canc{\`e}s, C., Chainais-Hillairet, C., and Krell, S. (2017a).
\newblock A nonlinear discrete duality finite volume scheme for
  convection-diffusion equations.
\newblock In {\em International Conference on Finite Volumes for Complex
  Applications}, pages 439--447. Springer.

\bibitem[Canc{\`e}s et~al., 2017b]{cances2017numerical}
Canc{\`e}s, C., Chainais-Hillairet, C., and Krell, S. (2017b).
\newblock Numerical analysis of a nonlinear free-energy diminishing discrete
  duality finite volume scheme for convection diffusion equations.
\newblock {\em arXiv preprint arXiv:1705.10558}.

\bibitem[Etchegaray, 2016]{etchegaray:tel-01533458}
Etchegaray, C. (2016).
\newblock {\em {Mathematical and numerical modelling of cell migration}}.
\newblock Theses, {Universit{\'e} Paris-Saclay}.

\bibitem[Etchegaray et~al., 2017a]{JMB2D}
Etchegaray, C., Meunier, N., and Voituriez, R. (2017a).
\newblock A {2D} deterministic model for cell crawling, a minimal multiscale
  approach.
\newblock Submitted.

\bibitem[Etchegaray et~al., 2017b]{etchegaray2017analysis}
Etchegaray, C., Meunier, N., and Voituriez, R. (2017b).
\newblock Analysis of a non-local and non-linear fokker-planck model for cell
  crawling migration.
\newblock {\em arXiv preprint arXiv:1701.06862}.

\bibitem[Joanny and Prost, 2009]{Joanny2009Active}
Joanny, J.-F. and Prost, J. (2009).
\newblock Active gels as a description of the actin-myosin cytoskeleton.
\newblock {\em HFSP J}, 3(2):94--104.

\bibitem[Kruse et~al., 2005]{Kruse2005Generic_theory}
Kruse, K., Joanny, J.~F., J{\"u}licher, F., Prost, J., and Sekimoto, K. (2005).
\newblock Generic theory of active polar gels: a paradigm for cytoskeletal
  dynamics.
\newblock {\em Eur Phys J E Soft Matter}, 16(1):5--16.

\bibitem[Maiuri et~al., 2015]{Maiuri2015Actin-flows-med}
Maiuri, P., Rupprecht, J.-F., Wieser, S., Ruprecht, V., B{\'e}nichou, O.,
  Carpi, N., Coppey, M., De~Beco, S., Gov, N., Heisenberg, C.-P., Lage~Crespo,
  C., Lautenschlaeger, F., Le~Berre, M., Lennon-Dumenil, A.-M., Raab, M.,
  Thiam, H.-R., Piel, M., Sixt, M., and Voituriez, R. (2015).
\newblock Actin flows mediate a universal coupling between cell speed and cell
  persistence.
\newblock {\em Cell}, 161(2):374--86.

\end{thebibliography}
\end{document}